\documentclass[smallextended]{svjour3}       
\smartqed  % flush right qed marks, e.g. at end of proof

\usepackage{graphicx}
\usepackage{mathptmx}      % use Times fonts if available on your TeX system
\usepackage{amsmath}
\usepackage{amsfonts}

\def  \P{\mathbb {P}}
\def  \R{\mathbb {R}}

\journalname{Journal of Theoretical Probability}

\begin{document}

%%%%%%%%%%%%%%%%%%%%%%%%%%%%%%%%%%%%%%%%%%%%%%%%%%
%%%%%%%%%%%%%%%%          TITLE         %%%%%%%%%%%%%%%%%%%%%%%%%%
%%%%%%%%%%%%%%%%%%%%%%%%%%%%%%%%%%%%%%%%%%%%%%%%%%

\title{A Note on Distribution Free Symmetrization Inequalities}
\titlerunning{A Note on Two Probabilistic Symmetrization Inequalities}

\author{  Zhao Dong \and Jiange Li \and Wenbo V. Li }

\institute
{     
      Z. Dong   \at
              Institute of Applied Mathematics, Academy of Mathematics and Systems Sciences,
Chinese Academy of Sciences, Beijing, 100190.   \\
              \email {dzhao@amt.ac.cn} 
  \and                              
        J. Li   \at
              Department of Mathematical Sciences, University of Delaware, Newark, DE, 19716.\\             
              \email{lijiange@math.udel.edu} 
\and
       W.V. Li   \at
              Department of Mathematical Sciences, University of Delaware, Newark, DE, 19716. \\             
              \email{wli@math.udel.edu}    
}

\date{Received: date / Accepted: date}

\maketitle

%%%%%%%%%%%%%%%%%%%%%%%%%%%%%%%%%%%%%%%%%%%%%%%%%%
%%%%%%%%%%%%%%%%%%          ABSTRACT         %%%%%%%%%%%%%%%%%%%%%
%%%%%%%%%%%%%%%%%%%%%%%%%%%%%%%%%%%%%%%%%%%%%%%%%%

\begin{abstract}
Let $X, Y$ be two independent identically distributed (i.i.d.) random variables taking values from a separable Banach space $(\mathcal{X}, \|\cdot\|)$. Given two measurable subsets $F, K\subseteq\cal{X}$, we established distribution free comparison inequalities between $\P(X\pm Y \in F)$ and $\P(X-Y\in K)$. These estimates are optimal for real random variables as well as when $\mathcal{X}=\R^d$ is equipped with the $\|\cdot\|_\infty$ norm. Our approach for both problems extends techniques developed by Schultze and Weizs\"acher (2007).


\end{abstract}

%%%%%%%%%%%%%%%%%%%%%%%%%%%%%%%%%%%%%%%%%%%%%%%%%%
%%%%%%%%%%%%%%%%%%         INTROCUDTION         %%%%%%%%%%%%%%%%%%%
%%%%%%%%%%%%%%%%%%%%%%%%%%%%%%%%%%%%%%%%%%%%%%%%%%

\section{Introduction}
\label{intro}
Symmetrization is one of the most basic and powerful tools in probability theory, particularly in the study of sums of random variables, see Ledoux and Talagrand (1991). Using symmetrization techniques, some important results about symmetric random variables can be extended to more general situations, for example L\'{e}vy-type inequalities, L\'{e}vy-It\^{o}-Nisio theorem, etc. This note is motivated by a recent paper of Schultze and Weizs\"acher, in which they proved that for arbitrary symmetric random walk in $\R$ with independent increment the probability of crossing a level at given time $n$ is $O(n^{-1/2})$. The following distribution free symmetrization inequality played an important role in removing the symmetry assumption. For i.i.d. real random variables $X, Y$, they proved
\begin{eqnarray}
\P(|X+Y|\leq 1)<2\cdot\P(|X-Y|\leq 1).
\end{eqnarray}
As mentioned by Schultze and Weizs\"acher,  the mere existence of a symmetrization constant for higher dimensional version of (1) follows from the estimate (29) of Mattner (1996). In this note, we extend the inequality to the general Banach space setting. In this context, our extension is useful in investigating analogous phenomena of random walks in Banach space, although we do not have a close study of this problem in the current paper. Moreover, our estimates are tight for real random variables as well as when $\mathcal{X}=\R^d$ is equipped with the $\|\cdot\|_\infty$ norm. The same approach is used to generalize a result of Alon and Yuster (1995): for all i.i.d. real random variables $X, Y$,
\begin{eqnarray}
\P(|X-Y|\leq b)<(2\lceil b/a \rceil-1)\cdot \P(|X-Y|\leq a).
\end{eqnarray} 
This answered a question of G. A. Margulis and Y. Peres. Moreover, Alon and Yuster showed the connection between optimal constants in such inequalities and kissing numbers, which have a long history of study; the kissing number in $\R^3$ was a subject of discussion between Isaac Newton and David Gregory in 1694. In addition, such estimates can be used to obtain moment inequalities involving certain classes of functions of $X+Y$ and $X-Y$. These and other applications will be studied in our subsequent paper.

Let $(\mathcal{X},\|\cdot\|)$ be a separable Banach space, and $F, K$ be two subsets of $\cal{X}$. We denote by $F\backslash K$ the set consisting of all elements in $F$ but not in $K$. Their sum is defined by
$$
F+K:=\{a+b: a\in F,~ b\in K\}.
$$ 
For $\rho>0$, we define the $\rho$-covering number of $F$ by $K$ in the following way
\begin{eqnarray}
N(F, K, \rho):=\inf \{|A|: A\subseteq \mathcal{X},  F\subseteq A+\rho K\}.
\end{eqnarray}
The usual definition of the diameter of $K$ is 
\begin{eqnarray}
d(K):=\sup_{x,y\in K}\|x-y\|,
\end{eqnarray}
and the inner radius is defined by
\begin{eqnarray}
r(K):=\sup\{r\geq 0: B(r) \subseteq K\},
\end{eqnarray}
where $B(r)$ is the closed ball centered at the origin with radius $r$. In general, we denote by $B(x, r)$ the closed ball centered at $x$ with radius $r$.
\begin{theorem}
Let $X, Y$ be i.i.d. $\mathcal{X}$-valued random variables and $F, K$ be two measurable subsets. If $K$ is symmetric and $r(K)>0$, we have
\begin{eqnarray}
\P(X +Y \in F)\leq N(F,  K,  \rho_K) \cdot\P(X-Y\in K),
\end{eqnarray}
where $\rho_K=r(K)/d(K)$. If $F$ is also symmetric, we have
\begin{eqnarray}
\P(X -Y \in F)\leq \left[N(F\backslash K,  K, \rho_K)+1\right]\cdot\P(X-Y\in K).
\end{eqnarray}

\end{theorem}

Let $B_1(r), B_2(r)\subseteq \R^d$ be two closed balls centered at the origin with radius $r$ under any two norms $\|\cdot\|_1$ and $\|\cdot\|_2$, respectively. Using Theorem 1, we have
\begin{corollary} 
For $a,b>0$ and i.i.d. $\mathbb{R}^d$-valued random variables $X, Y$, we have
\begin{eqnarray}
\P(\|X +Y\|_2\leq b)\leq N(B_2(b),  B_1(a), 1/2)\cdot\P(\|X -Y\|_1\leq a),
\end{eqnarray}
and
\begin{eqnarray}
\P(\|X -Y\|_2\leq b)\leq [N(B_2(b)\backslash B_1(a),  B_1(a), 1/2)+1]\cdot\P(\|X-Y\|_1\leq a).
\end{eqnarray}
\end{corollary}

Corollary 1 is tight for real random variables and the strict inequalities hold. The extension of (1) is given in the following theorem.
\begin{theorem} 
For $0<a/2<b$ and i.i.d. real random
variables $X, Y$, we have
\begin{eqnarray}
\P(|X+Y|\leq b)<\lceil2b/a\rceil\cdot \P(|X-Y|\leq a).
\end{eqnarray}
Moreover, the constant $\lceil2b/a\rceil$ can not be improved. When $0<b\leq a/2$, the inequality is still tight with $``\leq"$ in the middle.
\end{theorem}

It is not hard to see that Theorem 2 and the estimate (2) imply the following sharp inequalities.
\begin{corollary}
For $0<a/2<b$ and i.i.d. $\R^d$-valued random vectors $X, Y$ with independent entries, we have
\begin{eqnarray}
\P(\|X +Y\|_\infty\leq b)< \left(\lceil 2b/a\rceil\right)^d\cdot\P(\|X -Y\|_\infty\leq a).
\end{eqnarray}
For all $a, b>0$, we have
\begin{eqnarray}
\P(\|X -Y\|_\infty\leq b)< \left(2\lceil b/a\rceil-1\right)^d\cdot\P(\|X -Y\|_\infty\leq a).
\end{eqnarray}
\end{corollary}

%%%%%%%%%%%%%%%%%%%%%%%%%%%%%%%%%%%%%%%%%%%%%%
%%%%%%%%%%%%        PROOF OF MAIN THEOREM  1    %%%%%%%%%%%%%%%
%%%%%%%%%%%%%%%%%%%%%%%%%%%%%%%%%%%%%%%%%%%%%%

\section{ Proof of Theorem 1}
\label{sec:1}
The following lemma was proved by Schultze and Weizs\"acher, which shows how to derive two-variable inequalities from one-variable estimate.  We state the lemma in a form suitable for our purpose.
\begin {lemma}

Let $(\Omega,\mathcal {B})$ be a measurable space
and $f: \Omega\times\Omega\rightarrow\mathbb{R}$ be a
$\mathcal{B}\otimes\mathcal{B}$ measurable bounded symmetric
function. Let $\mathcal{P}$ be the set of all probability measures
on $\mathcal{B}$. Then the following statements are equivalent:

\begin{itemize}

\item For all $\mu\in\mathcal{P}$,
$$
\int_{\Omega\times\Omega}f(x,y)d\mu(x)d\mu(y)>0.
$$

\item For all $\mu\in\mathcal{P}$,
$$
\mu\left(\left\{x\in \Omega: \int_{\Omega}f(x,y)d\mu(y)>0\right\}\right)>0.
$$

\end{itemize}
\end{lemma}

\begin{proof}(Theorem 1):
We use $\mathcal{P}$ to denote the set of all probability measures on $\cal{X}$. Without confusion, we let $\rho:=\rho_K$ and $N:= N(F,  K,  \rho)$. Apparently, the theorem is true for $N=\infty$. In the following, we always assume $N$ is finite. In order to prove (6), 
we only need to show that for any constant $C>N$,
\begin{eqnarray}
\P(X + Y \in F)< C\cdot\P(X-Y\in K)
\end{eqnarray}
for all i.i.d. random variables $X, Y$. The inequality above can be rewritten as
\begin{eqnarray}
\int_{\cal{X}\times \cal{X}}\varphi(x,y)d\mu(x)d\mu(y)>0,
\end{eqnarray}
where 
$$
\varphi(x,y)=C\cdot1_{\{(x,y):  x-y\in K\}}-1_{\{(x,y):  x + y \in F\}},~~x,y\in \cal{X},
$$
and $\mu\in\mathcal{P}$ is induced by $X$.  Since $K$ is symmetric, we can see $\varphi(x,y)$ is symmetric and bounded.  By Lemma 1, it is equivalent to prove
\begin{eqnarray}
\mu\left(\left\{x\in \mathcal{X}: \int_{\mathcal{X}}\varphi(x,y)d\mu(y)>0\right\}\right)>0
\end{eqnarray}
for all $\mu\in\mathcal{P}$. Assume otherwise, then there exists some $\mu\in\mathcal{P}$ such that $\mu(S)=1$, where
\begin{eqnarray}
S &=&\left\{x\in \mathcal{X}: \int_{\mathcal{X}} \varphi(x,y)d\mu(y)\leq0\right\}\nonumber \\
&=&\left\{x\in \mathcal{X}: \mu\left(-x+F\right) \geq C\cdot\mu\left(x-K\right)\right\}.
\end{eqnarray}
Let's define
\begin{eqnarray}
\alpha=\sup_{x\in S}\mu\left(x-K \right).
\end{eqnarray}
Since $r(K)>0$ and $\mathcal{X}$ is separable, there exists a countable subset $S'\subseteq S$ such that $S\subseteq S'-K:= \cup_{x\in S'}(x-K)$, which implies $\alpha>0$. For $\epsilon>0$ small, we can pick $x^*\in S$ such that
\begin{eqnarray}
\mu(x^*-K)>\alpha-\epsilon.
\end{eqnarray}
By the definition of $N$, there exists a subset $\{x_i\}_{i=1}^{N}\subseteq \mathcal{X}$ such that
$$
F\subseteq\cup_{i=1}^{N}(x_i+\rho K).
$$
So we have
\begin{eqnarray}
-x^*+F\subseteq\cup_{i=1}^{N}(x_i-x^*+\rho K)=\cup_{i=1}^{N}(x_i-x^*-\rho K).
\end{eqnarray}
From (18), (16) and (19), we have
\begin{eqnarray}
C\cdot(\alpha-\epsilon)<C\cdot \mu(x^*-K)\leq\mu(-x^*+F)\leq N\cdot\sup_{x\in \mathcal{X}}\mu(x-\rho K).
\end{eqnarray}
Since $\mu(S)=1$, for any set  $x-\rho K$ with positive measure, there is
\begin{eqnarray}
x_0 \in (x-\rho K)\cap S.
\end{eqnarray}
Next we will show 
\begin{eqnarray}
x-\rho K\subseteq B(x_0, r(K))\subseteq x_0-K.
\end{eqnarray}
By (21), there exists $y_0\in K$ such that $x_0=x-\rho y_0$. For any $y\in K$,
$$
\|x-\rho y-x_0\| = \rho\|y_0-y\|\leq \rho\cdot d(K)=r(K),
$$
which implies the first part of (22). The second part follows from the assumption on $K$ and the definition of $r(K)$.  Combining (20)-(22), we have
$$
C\cdot(\alpha-\epsilon)<N\cdot\sup_{x\in\mathcal{X}}\mu(x-\rho K)\leq N\cdot\sup_{x\in S}\mu(x-K).
$$
Taking $\epsilon=\alpha\cdot(1-N/C)$, we have
\begin{eqnarray}
N\cdot \alpha=C\cdot(\alpha-\epsilon)<N\cdot\sup_{x\in S}\mu(x-K),
\end{eqnarray}
which contradicts the definition of $\alpha$ in (17). So we proved (6).\\

To prove (7), we only need to make a slight modification of the previous proof. Similar to (13), we need to prove that for any $C>N:=N(F\backslash K, K, \rho)+1$,
\begin{eqnarray}
\P(X -Y \in F)< C\cdot\P(X-Y\in K).
\end{eqnarray}
Instead of (16), we redefine
\begin{eqnarray}
S=\left\{x\in \mathcal{E}: \mu\left(x-F\right) \geq C\cdot\mu\left(x-K\right)\right\},
\end{eqnarray}
and $\alpha$ is defined in the same way as in (17). For $\epsilon>0$ small, we can pick $x^*\in S$ such that
\begin{eqnarray}
\mu(x^*-K)>\alpha-\epsilon.
\end{eqnarray} 
By the definition of $N$, there exists a subset $\{x_i\}_{i=1}^{N-1}\subseteq \mathcal{X}$ such that
$$
F\backslash K\subseteq\cup_{i=1}^{N-1}(x_i+\rho K).
$$
Hence
\begin{eqnarray}
x^*-F\subseteq (x^*-K) \cup \left(\cup_{i=1}^{N-1}(x^*-x_i-\rho K)\right).
\end{eqnarray}
From (26), (25) and (27), we have
\begin{eqnarray}
C\cdot(\alpha-\epsilon)&<&C\cdot \mu(x^*-K)\leq\mu(x^*- F)\\
&\leq& \mu(x^*-K) +(N-1)\cdot\sup_{x\in \mathcal{X}}\mu(x-\rho K).
\end{eqnarray}
Combining (28), (29), (21) and (22), we have
$$
C\cdot(\alpha-\epsilon)<\mu(x^*-K) +(N-1)\cdot\sup_{x\in \mathcal{X}}\mu(x-\rho K)\leq N\cdot\sup_{x\in S}\mu(x-K).
$$
Taking $\epsilon=\alpha\cdot(1-N/C)$, we get (23) again, which is in contradiction to the definition of $\alpha$. So we  proved (7).

\end{proof}

%%%%%%%%%%%%%%%%%%%%%%%%%%%%%%%%%%%%%%%%%%%%%
%%%%%%%%%%%%        PROOF OF THEOREM 2            %%%%%%%%%%%%%%%
%%%%%%%%%%%%%%%%%%%%%%%%%%%%%%%%%%%%%%%%%%%%%

\section{Proof of Theorem 2}

For $F=[-b, b]$, $K=[-a, a]$, we can see $\rho(K)=1/2$ and $N(F, K, 1/2)=\lceil 2b/a \rceil$.  In this case, Theorem 1 implies a slight weaker version of Theorem 2 without the strict inequality in the middle of (10). When $0<b\leq a/2$, the following trivial example shows that the equality can indeed happen. When $X, Y$ have the same distribution $\P(X=0)=1$, it is easy to see $\P(|X+Y|\leq b)=\P(|X-Y|\leq a)=1$. For $0<a/2<b$, we will extend the proof of (1) by Schultze and Weizs\"acher in the following section. 
\subsection{Generalization}
Without loss of generality, we assume $a=1$. By Lemma 1, we only need to prove the following claim.
\begin{claim}
Let $\mu$ be the probability measure on $\R$ induced by $X$, and $\mu_r(x)$ is defined by
$$
\mu_r(x):=\mu\left([x-r,x+r]\right).
$$
Then we have
\begin{eqnarray*}
\mu\left(\{x\in\R: \mu_b(-x)<\lceil 2b\rceil\cdot\mu_1(x)\}\right)>0.
\end{eqnarray*}
\end{claim}

\begin{proof}
If the claim is not true, there is some $\mu$ such that $\mu(S)=1$, where
\begin{eqnarray}
S=\left\{x\in \R: \mu_b(-x)\geq\lceil 2b\rceil\cdot\mu_1(x)\right\}.
\end{eqnarray}
Define 
$\alpha=\sup_{x\in S}\mu_1(x)$, which is positive. For $\epsilon>0$ small, we will show that there exists a sequence of disjoint intervals $\{I_k\}$ such that
\begin{eqnarray}
\mu(I_k)>\alpha-\lceil 2b\rceil^{2k}\epsilon.
\end{eqnarray}
For $M$ large enough, we have
$$
\mu\left(\cup_{k=0}^MI_k\right)>\sum_{k=0}^M(\alpha-\lceil 2b\rceil^{2k}\epsilon)>1,
$$
which is impossible. So the claim must be true.  Firstly, we can pick $x_0\in S$ such that $\mu_1(x_0)>\alpha-\epsilon$, and $I_0$ is defined as
\begin{eqnarray}
I_0=[x_0-1,x_0+1].
\end{eqnarray}
Since $x_0\in S$, we have
$$
\mu_b(-x_0)>\lceil 2b\rceil(\alpha-\epsilon).
$$
Without loss of generality, we assume $x_0\geq0$. It is easy to see that $[-x_0-b, -x_0+b]$ can be divided into $\lceil 2b\rceil$ disjoint intervals of the form
$$
[-x_0+b-1, -x_0+b], [-x_0+b-2, -x_0+b-1), \cdots, [-x_0-b, -x_0+b+1-\lceil 2b\rceil).
$$
Due to $\mu(S)=1$, the interval above with positive measure must have non-empty intersection with $S$. So it can be covered by $[y-1, y+1]$ for some $y\in S$. Then we can see that every interval above has measure at most $\alpha$, which implies
$$
\mu\left([-x_0-b, -x_0+b+1-\lceil 2b\rceil)\right)>\lceil 2b\rceil(\alpha-\epsilon)-(\lceil 2b\rceil-1)\alpha=\alpha-\lceil 2b\rceil\epsilon.
$$
For any $x_1\in [-x_0-b, -x_0+b+1-\lceil 2b\rceil)\cap S$, we have $\mu_1(x_1)>\alpha-\lceil 2b\rceil\epsilon$ and
\begin{eqnarray}
\mu_b(-x_1)> \lceil 2b\rceil(\alpha-\lceil 2b\rceil\epsilon).
\end{eqnarray}
When $b>1/2$, we always have
\begin{eqnarray}
-x_1+b> x_0+\lceil 2b\rceil-1>x_0+1.
\end{eqnarray}
For $1/2<b\leq1$, we can see
\begin{eqnarray}
x_0\geq-x_1-b>x_0+\lceil 2b\rceil-2b-1\geq x_0-1.
\end{eqnarray}
Combining (33)-(35), we have
\begin{eqnarray}
\mu((x_0+1, -x_1+b])&\geq&\mu_b(-x_1)-\mu_1(x_0)>\alpha-\lceil 2b\rceil^2\epsilon.
\end{eqnarray}
For $b>1$, we have
\begin{eqnarray}
-x_1-b-1+\lceil 2b\rceil>x_0+2(\lceil 2b\rceil-b-1)\geq x_0+1.
\end{eqnarray}
In this case, we also have
\begin{eqnarray}
\mu((-x_1-b-1+\lceil 2b\rceil, -x_1+b])\geq \alpha-\lceil 2b\rceil^2\epsilon.
\end{eqnarray}
Hence, we can define
\begin{eqnarray}
I_1=
\left \{\!\!\!
\begin{array}{ll}
(x_0+1, -x_1+b]~~~~~~~~~~~~~~~~1/2<b\leq 1,\\
(-x_1-b-1+\lceil 2b\rceil, -x_1+b]~~~~~~b>1.
\end{array}
\right.
\end{eqnarray}
Apparently, we have $I_0\cap I_1=\emptyset$. Proceeding recursively we can construct a sequence of disjoint intervals $\{I_k\}$ with properties as we mentioned before. So, the claim is true.
\end{proof}

\subsection{Example} 

In the following, we construct an example which shows that our estimate in Theorem 2 is sharp. Let $X,Y$ be independent random variables with the same distribution $\P(X=x_i)=(2n)^{-1}$, where
$$
x_i= \left \{\!\!\!
\begin{array}{ll}
~i(1+\epsilon)a~~~~~~~~~~~~i=1,2,\cdots,n,\\
~i(1+\epsilon)a-r~~~~~~i=0,-1,\cdots,-n+1,
\end{array}
\right.
$$
with $\epsilon>0$ small and $0< r\leq a(1+\epsilon)/2$. It is easy to see
\begin{eqnarray}
\P(|X-Y|\leq a)=\P(X=Y)=(2n)^{-1},
\end{eqnarray}
and
\begin{eqnarray}
\P(|X+Y|\leq1)=(2n)^{-1}\Big(\sum_{i\in I_1}+\sum_{i\in I_2}+\sum_{i\in I_3}\Big)\P(-x_i-1\leq X\leq-x_i+1),
\end{eqnarray}
where $\{I_1,I_2,I_3\}$ is a partition of the index set $\{i: -n+1\leq i\leq n\}$. The sets $I_1, I_2$ are defined by
\begin{eqnarray*}
I_1&=&\{i: - x_0+1\leq x_i\leq -x_{-n+1}-1\},\\
I_2&=&\{i: -x_n+1\leq x_i\leq-x_1-1\}.
\end{eqnarray*}
Elementary calculations show that 
\begin{eqnarray}
|I_1| &=& \lfloor
n-1-(1-r)(1+\epsilon)^{-1}a^{-1}\rfloor-\lceil(1+r)(1+\epsilon)^{-1}a^{-1}\rceil+1,\\
|I_2| &=&  \lfloor
n-(1+r)(1+\epsilon)^{-1}a^{-1}\rfloor-\lceil1+(1-r)(1+\epsilon)^{-1}a^{-1}\rceil+1.
\end{eqnarray}
For any $i\in I_1\cup I_2$, we have
\begin{eqnarray}
&&\P(-x_i-1\leq X\leq-x_i+1)=(2n)^{-1}\cdot|\{k:-x_i-1\leq
x_k\leq-x_i+1\}| \nonumber\\
&=&(2n)^{-1}\cdot\left(1+\lfloor
(1-r)(1+\epsilon)^{-1}a^{-1}\rfloor+\lfloor(1+r)(1+\epsilon)^{-1}a^{-1}\rfloor\right).
\end{eqnarray}
For any $i\in I_3$, we can see
\begin{eqnarray}
\P(-1-x_i\leq
X\leq1-x_i)=O(n^{-1}).
\end{eqnarray}
Combining (40)-(45), we have
\begin{eqnarray}
\lim_{n\rightarrow\infty}\frac{\P(|X+Y|\leq1)}{\P(|X-Y|\leq a)}=1+\lfloor
(1-r)(1+\epsilon)^{-1}a^{-1}\rfloor+\lfloor(1+r)(1+\epsilon)^{-1}a^{-1}\rfloor.
\end{eqnarray}
For all $a>0$, we will see that there are always appropriate $\epsilon, r$ such that the ratio above can achieve $\lceil2/a\rceil$.
\begin{enumerate}

\item When $k<1/a\leq k+1/2$, for some non-negative integer $k$, and $r>0$ small, we have
$$
k<(1-r)a^{-1}<k+1,~~k<(1+r)a^{-1}<k+1.
$$
For $\epsilon>0$ small, we have
$$
1+\lfloor
(1-r)(1+\epsilon)^{-1}a^{-1}\rfloor+\lfloor(1+r)(1+\epsilon)^{-1}a^{-1}\rfloor=2k+1=\lceil2/a\rceil.
$$

\item When $k+1/2<1/a\leq k+1$, and $r=a/2$, we have
$$
k<(1-r)a^{-1}<k+1<(1+r)a^{-1}<k+2.
$$
Then we can choose $\epsilon>0$ small such that
$$
1+\lfloor
(1-r)(1+\epsilon)^{-1}a^{-1}\rfloor+\lfloor(1+r)(1+\epsilon)^{-1}a^{-1}\rfloor=2k+2=\lceil2/a\rceil.
$$
\end{enumerate}

%%%%%%%%%%%%%%%%%%%%%%%%%%%%%%%%%%%%
%%%%%%%%%%%%              ACKNOWLEDGEMENT     %%%%%%%%%%
%%%%%%%%%%%%%%%%%%%%%%%%%%%%%%%%%%%%

\begin{acknowledgement}

We are grateful to Dr. Mokshay Madiman for his valuable suggestions and comments in the preparation of this note.

\end{acknowledgement}

%%%%%%%%%%%%%%%%%%%%%%%%%%%%%%%%%%%%
%%%%%%%%%%%%              REFERENCES          %%%%%%%%%%
%%%%%%%%%%%%%%%%%%%%%%%%%%%%%%%%%%%%

\end{document}